# Θεωρήσεις για τη φύση της Γεωμετρίας

**Ιωάννης Ρίζος** ioarizos@uth.gr
**Νικόλαος Γκρέκας** ngkrekas@uth.gr
Πανεπιστήμιο Θεσσαλίας, Τμήμα Μαθηματικών

ΠΕΡΙΛΗΨΗ

*Στην εργασία αυτή συζητάμε από ιστορική και φιλοσοφική άποψη μια εκδοχή για το νόημα των πέντε αιτημάτων της Ευκλείδειας Γεωμετρίας και κάνουμε σύντομη αναφορά στον φορμαλισμό του D. Hilbert. Εξετάζουμε διαχρονικά το ερώτημα «τι είναι η Γεωμετρία» μελετώντας την κατάτμηση των γεωμετριών αλλά και την ενοποίησή τους μέσα από το Πρόγραμμα Erlangen του F. Klein. Παρουσιάζουμε περιληπτικά τρεις εκ των εννέα γεωμετριών Cayley-Klein και κλείνουμε με κάποιες διδακτικές προοπτικές.*

**Λέξεις κλειδιά:** Γεωμετρία, αξίωμα, μετρική, ενοποίηση, Πρόγραμμα Erlangen, Διδακτική των Μαθηματικών, Ιστορία των Μαθηματικών

# Notions about the nature of Geometry

**Ioannis Rizos** ioarizos@uth.gr
**Nikolaos Gkrekas** ngkrekas@uth.gr
University of Thessaly, Department of Mathematics

ABSTRACT

*In this paper we discuss, from a historical and philosophical point of view, a variation of the meaning of the five postulates in Euclidean Geometry and we make a short reference to D. Hilbert's formalism. We examine, throughout the ages, the question «what is Geometry» by studying the segmentation and the unification of various geometries, been introduced by F. Klein's Erlangen Program. We concisely present three out of nine Cayley-Klein geometries and then we conclude with some didactic perspectives.*

**Εισαγωγή**

Η αλληλεπίδραση μεταξύ των μεταβολών στη θεώρηση της Γεωμετρίας και των εξελίξεων στις φυσικές θεωρίες, η οποία έγινε εντονότερη μετά την επιστημονική επανάσταση, έχει τις ρίζες της στους Πυθαγορείους και την Πλατωνική Ακαδημία: «Η φυσική πραγματικότητα για να είναι ορθολογική οφείλει να ενσαρκώνει μια γεωμετρική δομή». Η σκέψη αυτή συμπυκνώνεται στον *Τίμαιο* (38c, 53c και 54a αντίστοιχα) με τις κυκλικές κινήσεις των πλανητών, την τριγωνική μορφή της ύλης και τα πέντε κανονικά στερεά, και διαπνέει το σύνολο σχεδόν των *Στοιχείων* του Ευκλείδη. Το εντυπωσιακό είναι ότι αυτή η αλληλεπίδραση μεταξύ γεωμετρικής θεώρησης και φυσικής θεωρίας, τοποθετημένη βέβαια σε διαφορετική φιλοσοφική και επιστημολογική βάση, πέρασε στον Γαλιλαίο, τον Νεύτωνα και τον Einstein, φτάνοντας μέχρι τις μέρες μας: Πόσες είναι οι διαστάσεις του σύμπαντος; Ποιες γεωμετρικές δομές περιγράφουν καλύτερα τροποποιημένες θεωρίες βαρύτητας; Πού οφείλονται οι ρυτιδώσεις του χωροχρονικού συνεχούς; είναι μερικά μόνο από τα ερωτήματα που βρίσκονται στην αιχμή της επιστημονικής έρευνας. Σε κάποια από αυτά φιλοδοξούν να δώσουν απάντηση μοντέρνες θεωρίες, όπως η Μορφοκλασματική (Fractal) Γεωμετρία και η Γεωμετρία Finsler, η οποία αποτελεί ουσιαστικά γενίκευση της Γεωμετρίας Riemann χωρίς περιορισμό στην τετραγωνική μορφή.

Ο παραπάνω διάλογος μεταξύ Γεωμετρίας και Φυσικής, μας βοηθάει να κατανοήσουμε καλύτερα τη φύση των γεωμετρικών αξιωμάτων και μας υποδεικνύει την απουσία μαθηματικού μοντέλου a priori κατάλληλου να περιγράψει το σύνολο του φυσικού κόσμου. Στην εργασία μας επιχειρούμε να φωτίσουμε μέρη αυτού του διαλόγου, θέτοντας ως αφετηρία τη μελέτη της φύσης των Ευκλείδειων αιτημάτων. Παραθέτουμε εναλλακτικές ερμηνείες των αιτημάτων αυτών και κάνουμε νύξεις στην αξιωματική θεμελίωση του Hilbert. Μας απασχολεί η συμβολή του Felix Klein στην ταξινόμηση της Γεωμετρίας και οι ενοποιητικές δυνάμεις του Προγράμματος Erlangen. Σκιαγραφούμε τρεις επίπεδες Γεωμετρίες με κύριο άξονα τη μετρική των αποστάσεων και κλείνουμε συζητώντας ορισμένες διδακτικές προεκτάσεις.

**1. Φιλοσοφικό πλαίσιο και ιστορικό υπόβαθρο**

Με τη σημερινή οπτική μπορεί να θεωρηθεί πως υπάρχουν, μεταξύ άλλων, δύο βασικές αφετηρίες από τις οποίες έχει τη δυνατότητα να εκκινήσει κανείς προκειμένου να διατυπώσει αξιώματα/ αρχές για μια φυσική θεωρία.

Η πρώτη είναι *η παρατήρηση*. Επί παραδείγματι ο Kepler αξιοποίησε μια μακρά ακολουθία παρατηρησιακών δεδομένων για να διατυπώσει τους τρεις νόμους του για την κίνηση των πλανητών. Στην περίπτωση αυτή, η θεωρία εδραιώνεται από το πείραμα και τη συνακόλουθη μαθηματική επεξεργασία και τυποποίηση. Ωστόσο, αν μια τέτοιου είδους θεωρία θέλει να έχει πληροφοριακό περιεχόμενο, οφείλει να διατρέχει μονίμως τον κίνδυνο να διαψευστεί (βλ. σχετ. Popper 2002).

Η δεύτερη δυνατή αφετηρία είναι *οι φιλοσοφικές πεποιθήσεις* των επιστημόνων. Για παράδειγμα η παράδοση του «σώζειν τα φαινόμενα», θέτοντας αξιώματα για την κίνηση των πλανητών επί τέλειων κύκλων κ.λπ., έπαιξε καθοριστικό ρόλο στη γέννηση και εξέλιξη του κλάδου της μαθηματικής Αστρονομίας, από τις ομόκεντρες σφαίρες του Εύδοξου και του Κάλλιππου μέχρι και το *Mysterium Cosmographicum* του πρώιμου Kepler. Στην περίπτωση αυτή, η θεωρία δεν έχει ανάγκη νομιμοποίησης από την παρατήρηση ή το πείραμα. Η ισχύς της είναι άχρονη και αναλλοίωτη. Κάτι τέτοιο βέβαια για να συμβεί προϋποθέτει μια συγκεκριμένη κοινωνική/ πολιτιστική δομή, στην οποία να δεσπόζουν υπερβατικές αξίες και μεταφυσικές αντιλήψεις. Οι αντιλήψεις αυτές, οι οποίες μπορούν να συνοψιστούν στην έννοια του *πλατωνισμού*, κυριάρχησαν στη μαθηματική εκπαίδευση και έφτασαν, ως έναν βαθμό τουλάχιστον, μέχρι τις μέρες μας (Σπύρου 2009).

Όμως, πώς ακριβώς τίθενται τα αξιώματα στη Γεωμετρία; Όπως είναι γνωστό ο Ευκλείδης θεμελιώνει το πρώτο βιβλίο των *Στοιχείων* επάνω σε πέντε αξιώματα που αναφέρονται σε "σημεία", "ευθείες" και "επίπεδα". Υπάρχει άραγε κάποια σχέση ανάμεσα στις παραπάνω οντότητες και στην εμπειρική γνώση που προκύπτει από την παρατήρηση; Ο Poincaré (1952, σ. 49) θεωρεί πως «*δεν πειραματιζόμαστε επί ιδανικών ευθειών ή κύκλων. Πειράματα μπορούμε να πραγματοποιήσουμε μόνο επί υλικών αντικειμένων*». Η άποψη αυτή συμπλέει με το ρεύμα που αργότερα ονομάστηκε *ιντουισιονισμός* και δείχνει τη διάσταση μεταξύ της Γεωμετρίας που θεμελιώνεται σε a priori δοσμένα αξιώματα και μιας εξέλιξης κατευθυνόμενης από το πείραμα, η οποία επιβάλλει περιορισμούς και επιτρέπει αναθεωρητικές έρευνες «επί των υποθέσεων που κείνται στα θεμέλια της Γεωμετρίας», όπως υποδεικνύει και ο τίτλος της διατριβής του Riemann (2016).

Για τον ίδιο τον Ευκλείδη ωστόσο, τα αξιώματα της Γεωμετρίας δεν είναι υποθέσεις αλλά αυταπόδεικτες αλήθειες συμφυείς με τον χώρο που μας περιβάλλει. Για παράδειγμα το πρώτο αίτημα εξασφαλίζει την *ύπαρξη* και την *κατασκευή* ευθύγραμμου τμήματος από ένα σημείο σε ένα άλλο, ενώ η *μοναδικότητα* προκύπτει σε συνδυασμό με την Πρόταση Ι.4 (Heath 1956, p. 195), κάτι δηλαδή που συμφωνεί απόλυτα με την "καθημερινή εμπειρία" και την "κοινή λογική". Επιπλέον, αν ληφθεί υπόψη ότι ο Ευκλείδης ήταν «*τῇ προαιρέσει Πλατωνικός καὶ τῇ φιλοσοφίᾳ ταύτῃ οἰκεῖος*» (Πρόκλος, 68,20) και ότι τα *Στοιχεία* δημιουργήθηκαν σε μεγάλο βαθμό εντός της Πλατωνικής Ακαδημίας (Πρόκλος, 64,3-68,23), τότε δεν θα ήταν ορθό να αποκλειστεί εντελώς το ενδεχόμενο τα *Στοιχεία* να έλαβαν την τελική τους μορφή για να παίξουν *και έναν ρόλο εννοιολογικού πλαισίου* εντός του οποίου να μπορεί να μελετηθεί η Πλατωνική φυσική φιλοσοφία.[1] Μια τέτοια θεώρηση είναι συμβατή με το περιεχόμενο του Πλατωνικού *Τίμαιου*, όπως και με τη γνωστή ρήση του Γαλιλαίου ότι «*στο μεγάλο ανοιχτό βιβλίο του Σύμπαντος η Φιλοσοφία είναι γραμμένη στη γλώσσα των Μαθηματικών*» (*Il Saggiatore*, 6), αν και ο Γαλιλαίος κατανοούσε κάπως διαφορετικά τους όρους "ποσότητα" και "μαθηματική επιστήμη" (J. Klein 1998, σ. 51).

Έτσι λοιπόν, στο παραπάνω πλαίσιο, μπορούμε να θεωρήσουμε ότι το δεύτερο αίτημα εισάγει κατά μοναδικό τρόπο την *απειρία* στην Ευκλείδεια Γεωμετρία ενώ το τρίτο, με το αυθαίρετα μικρό και το αυθαίρετα μεγάλο "διάστημα", διασφαλίζει τη *συνέχεια* του χώρου (Heath 1956, p. 199).

Ο χαρακτήρας του τέταρτου αιτήματος είναι διαφορετικός σε σχέση με εκείνον των προηγουμένων, τα οποία είναι «*αιτήματα ύπαρξης*» (Bunt κ.ά. 1981, σ.166) ή απλά έχουν σκοπό «*να μας αναγγείλουν ότι επιτρέπονται οι κατασκευές με κανόνα και διαβήτη*» (Blanché 1971, σ. 29). Καταρχάς για τον Ευκλείδη αποδεκτές γωνίες είναι μόνο οι ευθύγραμμες, διαφορετικά το αντίστροφο του αιτήματος (δηλαδή ότι «όλες οι γωνίες που είναι ίσες με μία ορθή είναι ορθές»), όπως απέδειξε ο Πάππος (Πρόκλος, 189,12), δεν ισχύει πάντοτε. Άλλωστε σε μη ευθύγραμμες γωνίες, όπως λ.χ. στις κερατοειδείς, δεν ισχύει η ιδιότητα Αρχιμήδους-Ευδόξου. Ως εκ τούτου, *ο χώρος δεν μπορεί να καμπυλώνεται* και δύο ορθές γωνίες να παίρνουν τη μορφή του σχήματος 1:

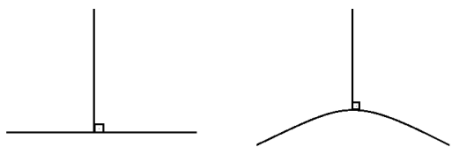

**Σχήμα 1:** Υποθετικά άνισες ορθές γωνίες

---

[1] Χωρίς να παραβλέπεται το γεγονός ότι τα *Στοιχεία*, παρά τις όποιες ατέλειές τους, αποτελούν μια εισαγωγή στα Μαθηματικά και είναι το πρώτο ολοκληρωμένο αξιωματικό σύστημα το οποίο περιλαμβάνει θεμελιώδη χαρακτηριστικά, όπως την ασυμμετρία κ.λπ.

Επιπλέον η ισότητα των δύο γωνιών (ως προς το μέγεθος και το σχήμα) ελέγχεται εάν τοποθετηθούν η μία επάνω στην άλλη. Η διαδικασία αυτή, γνωστή ως "εναπόθεση", αποτελεί ένα είδος κίνησης και προϋποθέτει ότι οι γωνίες –τα σχήματα εν γένει– παραμένουν αναλλοίωτα κατά τη μεταφορά. Η φύση όμως δεν μπορεί παρά να συμπεριφέρεται ομοιόμορφα, επομένως ο χώρος είναι εκ των ων ουκ άνευ *ομογενής*.[2] Κατά συνέπεια το τέταρτο αίτημα (ιδιαίτερα εάν ιδωθεί και σε συνδυασμό με την Πρόταση Ι.4) έχει σκοπό να εξασφαλίσει αυτή την ομογένεια. Συνοπτικά θα λέγαμε ότι το τέταρτο αίτημα παρέχει ένα *εγγενές πρότυπο μέτρησης γωνιών*, καθώς οι ορθές γωνίες ορίζονται γεωμετρικά και οι άλλες γωνίες έχουν πλέον τη δυνατότητα να συγκριθούν με αυτές (Greenberg 2008, p. 20).

Ο κατάλογος των αιτημάτων του βιβλίου Ι των *Στοιχείων* κλείνει με το πέμπτο και γνωστότερο "αίτημα των παραλλήλων" που, αν συγκριθεί το ύφος του με εκείνο των προηγουμένων, φαίνεται να το πρόσθεσε ο ίδιος ο Ευκλείδης προκειμένου να στηρίξει το οικοδόμημά του (Torretti 1984, σ. 43). Αίτημα ισχυρού διαισθητικού χαρακτήρα, με το οποίο, όπως και με το τέταρτο, γίνεται παραδεκτό ότι κάποια γεωμετρικά σχήματα έχουν ορισμένες ειδικές ιδιότητες. Γιατί όμως ο Ευκλείδης σταμάτησε εδώ, ενώ είναι πρόδηλο ότι το πλήθος των αιτημάτων του είναι περιορισμένο για να μπορέσει να θεμελιώσει αξιωματικά τη Γεωμετρία;

Η επικρατούσα –προφανής ομολογουμένως– εξήγηση είναι ότι βασίστηκε υπερβολικά στις αισθήσεις και θεώρησε αυτονόητες κάποιες ιδιότητες, με την έννοια ότι αυτές εκπορεύονται από ένα γνωστικά λειτουργικό σχήμα. Το κλασικότερο παράδειγμα είναι η απόδειξη της Πρότασης Ι.1 (*Ἐπὶ τῆς δοθείσης εὐθείας πεπερασμένης τρίγωνον ἰσόπλευρον συστήσασθαι*) όπου δεν προκύπτει από πουθενά ότι οι δύο κύκλοι που φέρνουμε από τα άκρα του ευθυγράμμου τμήματος τέμνονται. Άλλωστε σε ένα καρτεσιανό σύστημα με συντεταγμένες ρητούς αριθμούς, οι δύο κύκλοι ενδέχεται να μην έχουν κοινά σημεία. Για να θεωρηθεί λοιπόν η απόδειξη πλήρης, πρέπει να εισαχθεί ένα επιπλέον αξίωμα, παραλλαγή του αξιώματος του Pasch.

Δεν μπορεί να αποκλειστεί ωστόσο και μια φιλοσοφική προσέγγιση στο ερώτημα, εάν μάλιστα ληφθεί υπόψη ο τρόπος με τον οποίο οι αρχαίοι Έλληνες έβλεπαν τα αιτήματα. Συγκεκριμένα, ενδέχεται ο Ευκλείδης να παρέλειψε εκούσια κάποιες προτάσεις οι οποίες θα μπορούσαν να παίξουν

---

[2] Σήμερα λέμε ότι ο χώρος είναι *ομογενής* όταν σε ένα αδρανειακό σύστημα αναφοράς όλα τα σημεία του είναι ισοδύναμα. Έτσι, η μεταφορά ενός κλειστού συστήματος στον χώρο χωρίς αλλαγή των σχετικών θέσεων και ταχυτήτων των σωμάτων, αφήνει αναλλοίωτες τις μηχανικές ιδιότητες του συστήματος

τον ρόλο αιτημάτων διότι, προφανώς, είχε την πεποίθηση ότι *δεν χρειάζονται πολλά αιτήματα*.[3] Αντίθετα, οι ορισμοί του βιβλίου Ι είναι αρκετοί –23 συνολικά– και μάλλον περισσότεροι απ' όσους χρειάζονται, χωρίς όμως να ισχυρίζεται κανείς πως είναι απόλυτα διαφωτιστικοί. Για παράδειγμα το σημείο ορίζεται ως «*οὗ μέρος οὐθέν*», δηλαδή ως κάτι που δεν έχει διαστάσεις, κάτι άυλο! Έτσι, αν οι βασικές ιδιότητες του χώρου (συνέχεια, ομογένεια κ.λπ.) οι οποίες εισάγονται με τα πέντε αιτήματα μπορούσαν να οριστούν κατάλληλα, τότε όλα τα αξιώματα θα μετατρέπονταν σε αποδείξιμες προτάσεις. Η τάση αυτή θυμίζει έντονα τη σύγχρονη σχολή σκέψης του *λογικισμού*, σύμφωνα με την οποία τα αξιώματα και τα θεωρήματα μπορούν να αναχθούν στην καθαρή λογική. Ο Αριστοτέλης (*Τοπικά*, 158b24-159a2), χρησιμοποιώντας ως παράδειγμα τον ορισμό της αναλογίας, περιγράφει τη δυναμική σχέση μεταξύ αξιωμάτων και "καλών ορισμών". Μέσω αυτής της σχέσης μπορεί να υποστηριχθεί ότι ο βαθμός τελειότητας της θεωρίας λόγων μεγεθών του βιβλίου VII των *Στοιχείων*, είναι ανώτερος από τον βαθμό τελειότητας της Ευκλείδειας Γεωμετρίας του βιβλίου Ι (Νεγρεπόντης & Φαρμάκη 2019, σ. 86).

Εντούτοις η θεώρηση του Ευκλείδη δεν είναι η μόνη δυνατή αναφορικά με το πώς τίθενται τα αξιώματα. Από τις αρχές του 16ου αιώνα είχε ήδη αρχίσει να ασκείται κριτική –γόνιμη ως επί το πλείστον– στα *Στοιχεία*, η οποία σε συνδυασμό με τις προσπάθειες απόδειξης του πέμπτου αιτήματος, οδήγησε στην ανακάλυψη των μη Ευκλείδειων γεωμετριών από τους Gauss, Lobachevsky και Bolyai. Τότε μόνο «*έγινε σαφής η διάκριση μεταξύ της γεωμετρίας ως τυπικού αξιωματικού συστήματος και της γεωμετρίας ως πλαισίου μελέτης του φυσικού σύμπαντος*» (Davis 2001, σ. 93).

Η πρώτη καθαρά συμπερασματική κατασκευή της Γεωμετρίας έγινε από τον Hilbert, ο οποίος κινούμενος μέσα στο φιλοσοφικό πλαίσιο του *φορμαλισμού*, την οδήγησε σε «διαζύγιο από τη διαίσθηση», απεκδύοντας ουσιαστικά τους γεωμετρικούς όρους από το "παραδοσιακό" εννοιολογικό τους περιεχόμενο. Για παράδειγμα, η γραμμή έπαψε να είναι «*μῆκος ἀπλατές*» και έγινε, όπως το σημείο κ.τ.λ. *απροσδιόριστη έννοια*. Η διαίσθηση μολαταύτα, παρά το ότι εξαλείφεται από μια θεωρία –εν προκειμένω την Ευκλείδεια Γεωμετρία– δεν παύει να είναι αυτή η οποία υπαγόρευσε την εξέλιξη των σταδίων για την αξιωματικοποίηση της

---

[3] Χαρακτηριστικό είναι το γεγονός ότι ο Ευκλείδης φαίνεται να ανέβαλλε όσο μπορούσε περισσότερο τη χρησιμοποίηση του πέμπτου αιτήματος, το οποίο αξιοποιείται για πρώτη φορά στα *Στοιχεία* στην Πρόταση Ι.29, παρότι θα απλοποιούσε και κάποιες από τις προηγούμενες αποδείξεις.

θεωρίας (Blanche 1971) παίζοντας, θα μπορούσε να πει κανείς, τον ρόλο της "σκαλωσιάς".

Ο Hilbert έθεσε 20 αξιώματα σε 5 ομάδες (σύμπτωσης, διάταξης, ισότητας, συνέχειας, παραλληλίας), πιστεύοντας ότι με αυτή την αυστηρότητα θα θεμελίωνε ακλόνητα το γεωμετρικό του οικοδόμημα και θα απάλλασσε οριστικά τα Μαθηματικά από κάθε εσωτερική αντίφαση ή αντινομία. Σύμφωνα με τον Hilbert (1899), λοιπόν, ένα αξιωματικό σύστημα πρέπει να είναι: α) *Συνεπές*, δηλαδή να μην οδηγεί σε αντιφάσεις. Αυτό σημαίνει ότι είναι αδύνατον στο πλαίσιο του αξιωματικού συστήματος να αποδειχθεί ότι μια πρόταση είναι ταυτόχρονα αληθής και ψευδής β) *Ανεξάρτητο*, δηλαδή να μην είναι δυνατόν να αποδειχθεί ένα αξίωμα από τα υπόλοιπα γ) *Πλήρες*, δηλαδή για κάθε πρόταση να μπορεί να αποδειχθεί η αλήθεια ή το ψεύδος της. Αυτό ισοδύναμα σημαίνει ότι δεν μπορούμε να προσθέσουμε στο σύστημα ένα άλλο αξίωμα ανεξάρτητο από τα ήδη υπάρχοντα.

Το 1931 ο Gödel ανέτρεψε τις παραπάνω προϋποθέσεις, γκρεμίζοντας οριστικά και αμετάκλητα το «όνειρο του Hilbert».[4] Απέδειξε ότι σε οποιαδήποτε αξιωματική θεωρία που περιέχει την Peano αριθμητική: α) δεν είναι δυνατόν να αποδειχθεί μέσα στη θεωρία ότι αυτή δεν οδηγεί σε αντιφάσεις και β) μπορεί να διατυπωθεί πρόταση που ούτε η ίδια ούτε η άρνησή της να μπορούν να αποδειχθούν. Κάπως έτσι έκλεισε ένας κύκλος ζυμώσεων 2.200 ετών, ο οποίος –ουσιαστικά– είχε στο κέντρο του τις διαμάχες των μελών της επιστημονικής κοινότητας για την εποπτεία του χώρου και τη Γεωμετρία που τον περιγράφει.

## 2. Τι είναι, τελικά, η Γεωμετρία;

Εάν θελήσουμε να εξετάσουμε διαχρονικά το ερώτημα «τι είναι η Γεωμετρία;», θα διαπιστώσουμε πως είναι δύσκολο να δώσουμε μία μοναδική απάντηση. Διότι από τις πρακτικές εφαρμογές που μας ιστορεί ο Ηρόδοτος (*Ἱστορίαι*, ΙΙ,109) για τη «μέτρηση της γης» στην αρχαία Αίγυπτο, η Γεωμετρία πέρασε σε μια πρώτη φάση θεωρητικοποίησης με τον Θαλή και αξιωματικοποίησης –με τον τρόπο που είδαμε στην προηγούμενη παράγραφο– με τον Ευκλείδη. Έτσι, για μεγάλο χρονικό διάστημα υπήρχε η πεποίθηση πως η Γεωμετρία είναι «η μελέτη των επίπεδων και στερεών σχημάτων», μέχρι που τον 17ο αιώνα ο Descartes με τη *Γεωμετρία* και ο Fermat με την *Εισαγωγή στους επίπεδους και στερεούς τόπους*, έδειξαν πως η Γεωμετρία μπορεί να βασιστεί στην έννοια του αριθμού και να μετατραπεί σε Ανάλυση με την εισαγωγή ενός συστήματος συντεταγμένων στο επίπεδο. Έως τα τέλη του 18ου αιώνα και τον καταλυτικό ρόλο που διαδραμάτισε ο Kant, είχε καθιερωθεί η άποψη ότι η Γεωμετρία αποτελεί «τη μελέτη του χώρου», του οποίου η γνώση δεν μπορεί να είναι εμπειρική αλλά συνθετική a priori, συνέπεια του τρόπου δόμησης του νου.

---

[4] Η ειρωνεία είναι ότι την ίδια περίπου εποχή, στις 8 Σεπτεμβρίου 1930, ο Hilbert μιλώντας στη Γερμανική Εταιρεία Φυσικών Επιστημόνων (Gesellschaft der Deutschen Naturforscher und Ärzte) είχε δηλώσει: «*Δεν υπάρχει άγνοια* (ignorabimus) *για εμάς* [τους μαθηματικούς] *και κατά τη γνώμη μου καθόλου για την επιστήμη. Αντί για την ανόητη άγνοια, το σύνθημά μας οφείλει να είναι: Πρέπει να μάθουμε – θα μάθουμε* (Wir müssen wissen – wir werden wissen)». Η τελευταία φράση είναι χαραγμένη επάνω στον τάφο του, στο Göttingen.

Τον 19ο αιώνα η κατάσταση στη Γεωμετρία είχε περίπου ως εξής: Από τη μια μεριά υπήρχαν οι συνθετικοί γεωμέτρες, βασιζόμενοι σε καθαρά γεωμετρικές μεθόδους και σημείο εκκίνησης τα αξιώματα, όπως ο Poncelet, ο Steiner και ο von Staudt, και από την άλλη οι αναλυτικοί γεωμέτρες, στηριζόμενοι σε αλγεβρικές πρακτικές, όπως ο Plücker, ο Cayley και ο Grassmann. Υπήρχε η Παραστατική και η Προβολική Γεωμετρία, είχε διατυπωθεί η «αρχή του δυϊσμού», ενώ το αρχαίο ερώτημα για την ανεξαρτησία του πέμπτου αιτήματος είχε απαντηθεί από τους Gauss, Lobachevsky και Bolyai κατακερματίζοντας ακόμη περισσότερο τη Γεωμετρία (Millman 1977). Η μόνη προσπάθεια ενοποίησης έως τότε ήταν η διάλεξη (Habilitationsschrift) του Riemann το 1854 στο Πανεπιστήμιο του Göttingen (Riemann 2016), στην οποία διατυπώθηκε για πρώτη φορά η έννοια του «μαθηματικού χώρου» ως συνεχούς συνόλου ομογενών αντικειμένων ή φαινομένων. Η διάλεξη όμως εκείνη θεωρήθηκε, τότε, από πολλούς μαθηματικούς αρκετά αφηρημένη και δύσκολη, με συνέπεια να μην της δοθεί η προσοχή που της άξιζε. Έχοντας αυτά υπόψη, μπορεί πλέον να αποτιμηθεί η σημαντικότητα του Erlangen Program του Γερμανού μαθηματικού Felix Klein (1849-1925).

Ο Klein ήταν βοηθός του Julius Plücker στο Πανεπιστήμιο της Βόννης, από τον οποίο όμως δεν κληρονόμησε τον ενθουσιασμό για την Αναλυτική Γεωμετρία. Οι εργασίες του ακολούθησαν άλλον δρόμο – αυτόν που στόχευε να φέρει στοιχεία ενότητας στην ποικιλία των νέων αποτελεσμάτων της έρευνας στη Γεωμετρία (Boyer 1968, p. 591). Η διαφορετική αυτή οπτική ενδέχεται να οφείλεται εν μέρει στη διαμονή του Klein για ένα διάστημα στο Παρίσι, το 1870, όπου η Θεωρία Ομάδων είχε ήδη εξελιχθεί σε έναν αυτόνομο κλάδο της Άλγεβρας. Εκεί γνωρίστηκε με τον Νορβηγό μα-θηματικό Sophus Lie και συνεργάστηκε μαζί του, εφόσον και οι δύο συνειδητοποίησαν τη σημαντικότητα της έννοιας της ομάδας[5] στα Μαθηματικά.

Το 1872 ο Klein διορίζεται καθηγητής στο Πανεπιστήμιο του Erlangen στη Νυρεμβέργη. Εκεί δίνει μια διάλεξη εγκαινιάζοντας ουσιαστικά ένα ερευνητικό πρόγραμμα, αναδεικνύοντας την αξία των μετασχηματισμών και εξηγώντας το πώς η έννοια της ομάδας θα μπορούσε να αξιοποιηθεί ως ένα κατάλληλο μέσο ταξινόμησης των διάφορων γεωμετριών, οι οποίες είχαν κάνει την εμφάνισή τους τον 19ο αιώνα. Το πρόγραμμα αυτό, το οποίο έμεινε γνωστό ως «Erlangen Program», περιγράφει τη Γεωμετρία ως «τη μελέτη των ιδιοτήτων των σχημάτων οι οποίες παραμένουν αναλλοίωτες κάτω από μια συγκεκριμένη ομάδα μετασχηματισμών». Έτσι οποιαδήποτε ταξινόμηση ομάδων μετασχηματισμών, μετατρέπεται σε κωδικοποίηση γεωμετριών. Η σχέση μεταξύ των σχημάτων είναι *σχέση ισοδυναμίας*, δηλαδή αυτοπαθής, συμμετρική και μεταβατική. Άρα τα σχήματα ταξινομούνται σε κλάσεις ισοδυναμίας, με μοναδική διαφορά τη θέση τους σε κάθε γεωμετρία. Μαθηματικά αυτό μεταφράζεται ως εξής:

Έστω $G$ μια ομάδα και $X$ ένα μη κενό σύνολο. Θεωρούμε την απεικόνιση $\delta: G \times X \to X$ με $\delta(g, x) = g \cdot x$ τέτοια ώστε:
  i. $\delta(e, x) = x$, όπου $e$ το ουδέτερο στοιχείο της ομάδας $G$
  ii. $\delta(g_1, \delta(g_2, x)) = \delta(g_1 g_2, x)$

Αν ισχύουν τα παραπάνω, τότε η απεικόνιση $\delta$ καλείται *δράση* της ομάδας $G$ επί του συνόλου $X$. Το ζεύγος $(G, X)$ καλείται *Γεωμετρία κατά Klein*.

Για παράδειγμα, η Ευκλείδεια Γεωμετρία (στις δύο διαστάσεις) είναι η μελέτη εκείνων των ιδιοτήτων των σχημάτων, συμπεριλαμβανομένου του μήκους και του εμβαδού, που παραμένουν αναλλοίωτες κάτω από την ομάδα των μετασχηματισμών η οποία αποτελείται από *μεταφορές* και *στροφές* στο επίπεδο – ομάδα ουσιαστικά ισοδύναμη με το τέταρτο αίτημα του Ευκλείδη, τουλάχιστον με το περιεχόμενο που του αποδώσαμε στην προηγούμενη παράγραφο.

---

[5] Για την αυτοτέλεια του κειμένου, παραθέτουμε εδώ τον ορισμό της ομάδας: Έστω ένα μη κενό σύνολο $G$ εφοδιασμένο με μία πράξη $*$. Η δομή $(G, *)$ λέγεται *ομάδα* αν ισχύουν οι ιδιότητες: i) Αν $x, y \in G$, τότε $x * y \in G$ (κλειστότητα) ii) $x * (y * z) = (x * y) * z$ (προσεταιριστικότητα) iii) Υπάρχει $e \in G$ τέτοιο ώστε $x * e = e * x = x$, για κάθε $x \in G$ (ύπαρξη ουδέτερου στοιχείου) iv) Για κάθε $x \in G$ υπάρχει $x' \in G$ τέτοιο ώστε $x * x' = x' * x = e$ (ύπαρξη αντίθετου στοιχείου). Αν επιπλέον ισχύει ότι $x * y = y * x$ για όλα τα $x, y \in G$, τότε η ομάδα λέγεται *αντιμεταθετική ή αβελιανή*.

Ο Klein, στο παραπάνω πλαίσιο και στηριζόμενος στην προηγούμενη δουλειά του Cayley (1859), κατέταξε τις γεωμετρίες σε εννέα ομάδες, γνωστές σήμερα ως «Γεωμετρίες Cayley-Klein».[6] Επιλέγοντας έναν από τους τρεις (ελλειπτικό, παραβολικό, υπερβολικό) τρόπους μέτρησης του μήκους μιας ευθείας και έναν από τους τρεις τρόπους μέτρησης των γωνιών, μπορούμε να πάρουμε εννέα επίπεδες γεωμετρίες, όπως φαίνεται στον ακόλουθο πίνακα (Yaglom 1979, p. 218):

**Πίνακας 1** Εννέα Γεωμετρίες Cayley-Klein στο επίπεδο

| | | Μέτρο μήκους | | |
|---|---|---|---|---|
| | | Ελλειπτικό | Παραβολικό | Υπερβολικό |
| Μέτρο γωνιών | Ελλειπτικό | Ελλειπτική γεωμετρία | Ευκλείδεια γεωμετρία | Υπερβολική γεωμετρία |
| | Παραβολικό | συν-Ευκλείδεια γεωμετρία | Γεωμετρία Galileo | συν-Minkowski γεωμετρία |
| | Υπερβολικό | συν-υπερβολική γεωμετρία | Γεωμετρία Minkowski | Διπλή υπερβολική γεωμετρία |

Η σημασία του Προγράμματος Erlangen και ο ρόλος της Γεωμετρίας ως ενοποιητικής δύναμης όλων των μαθηματικών περιοχών, τονίστηκε από τον Κωνσταντίνο Καραθεοδωρή σε εργασία του αφιερωμένη στα $70^α$ γενέθλια του Felix Klein: «*Στο Πρόγραμμα Erlangen εμφανίστηκε για πρώτη φορά μια τάση η οποία συνίστατο στη σύνδεση ανάμεσα σε απομακρυσμένες περιοχές* [των Μαθηματικών], *δημιουργώντας έτσι νέες καρποφόρες δυνατότητες έρευνας και συμβάλλοντας στην υπέρβαση των κινδύνων κατακερματισμού της επιστήμης*» (Carathéodory 1919).[7]

Από τις γεωμετρίες Cayley-Klein, για τους σκοπούς της μελέτης μας, θα παρουσιάσουμε με συντομία τρεις στην επόμενη παράγραφο.

---

[6] Για έναν ιστορικό απολογισμό της ανακάλυψης αυτών των γεωμετριών βλ. Rosenfeld, B. A. (1988). *A History of Non-Euclidean Geometry*. New York: Springer-Verlag, το κλασικό Lewis, F. P. (1920). History of the parallel postulate. *American Mathematical Monthly*, 27(1), 16-23, καθώς και την άποψη του ίδιου του Klein στο Klein, F. (1928). *Vorlesungen über Nicht-Euklidische Geometrie*. Berlin: Springer.

[7] Μια πιο σύγχρονη αξιολόγηση της θέσης του Προγράμματος Erlangen στην Ιστορία των Μαθηματικών γίνεται στο Hawkins, T. (1984). The Erlanger Programm of Felix Klein: Reflections on Its Place in the History of Mathematics. *Historia Mathematica*, 11, 442-470. Για τη μετάβαση στη Γεωμετρία κατά E. Cartan βλ. Sharpe, R. W. (1997). *Differential Geometry: Cartan's Generalization of Klein's Erlangen Program*. New York: Springer.

**3. Τρεις επίπεδες γεωμετρίες**

Ακολουθώντας την άποψη του Klein για τη Γεωμετρία και αποφεύγοντας να δεχτούμε ως αξίωμα κάτι το οποίο μπορεί να προκύψει από τη μελέτη των φυσικών φαινομένων (π.χ. από τη μελέτη της κίνησης ρευματοφόρου αγωγού σε μαγνητικό πεδίο), παρουσιάζουμε συνοπτικά εδώ ορισμένα βασικά στοιχεία τριών γεωμετριών στις δύο διαστάσεις: της Ευκλείδειας Γεωμετρίας, της Γεωμετρίας Galileo και της Γεωμετρίας Minkowski.

Στην αναλυτική επίπεδη γεωμετρία, οι βασικές έννοιες είναι οι μετρικές των αποστάσεων και των γωνιών. Κάθε τι άλλο εξαρτάται από αυτές. Το απαραίτητο εργαλείο γι' αυτές τις μετρικές είναι ένα *εσωτερικό γινόμενο* στον $\mathbb{R}^2$. Το σύνηθες Ευκλείδειο εσωτερικό γινόμενο δύο διανυσμάτων $a = (x_1, y_1)$ και $b = (x_2, y_2)$ είναι:

$$<,>: \mathbb{R}^2 \times \mathbb{R}^2 \to \mathbb{R}: (a,b) \mapsto <a,b> := x_1 x_2 + y_1 y_2$$

Από το εσωτερικό γινόμενο λαμβάνουμε το μήκος ή *Ευκλείδεια νόρμα* ενός διανύσματος $u = (x, y)$ ως

$$\|u\| = \sqrt{<u,u>} = \sqrt{x^2 + y^2}.$$

Τώρα είμαστε σε θέση να κάνουμε τις βασικές μετρήσεις. Για δύο σημεία, ή καλύτερα για δύο διανύσματα $a = (x_1, y_1)$ και $b = (x_2, y_2)$ του $\mathbb{R}^2$, η απόστασή τους ορίζεται ως

$$\|a - b\| = \sqrt{(x_1 - x_2)^2 + (y_1 - y_2)^2}$$

ενώ η γωνία $\varphi$ που σχηματίζουν, δίνεται από τον τύπο

$$\varphi = arccos \frac{a \cdot b}{\|a\|\|b\|}.$$

Είναι θεμελιώδης γεωμετρική αρχή τα μήκη και οι γωνίες να μην αλλάζουν όταν το επίπεδο υφίσταται στροφή. Η εν λόγω στροφή υπό γωνία $\theta$ δίνεται από τον πίνακα

$$\begin{pmatrix} cos\theta & -sin\theta \\ sin\theta & cos\theta \end{pmatrix}.$$

Αποδεικνύεται ότι ο παραπάνω πίνακας διατηρεί την ποσότητα $x^2 + y^2$ *αναλλοίωτη*, δηλαδή αν $x^2 + y^2 = m$ και $\begin{pmatrix} x' \\ y' \end{pmatrix} = \begin{pmatrix} cos\theta & -sin\theta \\ sin\theta & cos\theta \end{pmatrix} \begin{pmatrix} x \\ y \end{pmatrix}$, τότε και $(x')^2 + (y')^2 = m$.[8]

Για τη μελέτη της «Γεωμετρίας της αρχής της σχετικότητας του Galileo» ή πιο απλά της «Γεωμετρίας Galileo», θέτουμε ως αφετηρία την εξής αρχή:

«Δύο παρατηρητές κινούμενοι ομαλά, ο ένας ως προς τον άλλον, θα πρέπει να αντιλαμβάνονται τους φυσικούς νόμους με τον ίδιο ακριβώς τρόπο. Κανένας παρατηρητής δεν μπορεί να διακρίνει την απόλυτη ηρεμία από την απόλυτη κίνηση, σε όποιον φυσικό νόμο κι αν προσφύγει. Γι' αυτό τον λόγο δεν υπάρχει *απόλυτη* κίνηση, αλλά μόνο κίνηση *σχετική* με κάποιον παρατηρητή».[9]

Μαθηματικοποιώντας την παραπάνω αρχή, θεωρούμε δύο αδρανειακά συστήματα αναφοράς (Α.Σ.Α.) τα οποία μελετούν την κίνηση σώματος μάζας $m$, με το δεύτερο σύστημα να κινείται ως προς το πρώτο με σταθερή ταχύτητα $u$. Για τις δυνάμεις που ασκούνται στο σώμα, μετρούμενες από τα δύο Α.Σ.Α., ισχύει (στον απόλυτο *χώρο* και *χρόνο* της Νευτώνειας Μηχανικής) ότι $\Sigma F = \Sigma F'$ δηλαδή

$$m \frac{d^2 r}{dt^2} = m \frac{d^2 r'}{dt'^2}$$

όπου $(r, t)$ και $(r', t')$ η θέση και ο χρόνος του σώματος στα δύο συστήματα αντίστοιχα. Σε περίπτωση μονοδιάστατης κίνησης (π.χ. στον άξονα των $x$) και με σχετική ταχύτητα $u$ του Σ' ως προς το Σ, που δεν είναι συγκρίσιμη με την ταχύτητα του φωτός, η παραπάνω ισότητα των δυνάμεων ισχύει αν και μόνο αν οι συντεταγμένες συνδέονται μεταξύ τους με τους λεγόμενους «μετασχηματισμούς Galileo»:

$$x' = x - ut, \quad y' = y, \quad z' = z, \quad t' = t.$$

---

[8] Περισσότερα για τους μετασχηματισμούς στο επίπεδο παραπέμπουμε στο απλό και διδακτικό βιβλίο Pettofrezzo, A. J. (1978). *Matrices and Transformations*. New York: Dover.
[9] Ο Γαλιλαίος διατυπώνει την αρχή αυτή χρησιμοποιώντας ένα λογοτεχνικά θαυμάσιο ύφος στο έργο του *Dialogo sopra i due massimi sistemi del mondo* (1632).

Αν τώρα περιοριστούμε μόνο στα $x$ και $t$ και συμβολίσουμε με $y$ τη χωρική συντεταγμένη και με $x$ τη χρονική ενός σημείου (σημειακής μάζας) κατά τη μονοδιάστατη κίνησή του, τότε οι παραπάνω μετασχηματισμοί γίνονται:

$$y' = y - ux, \quad x' = x.$$

Στη γλώσσα της Γραμμικής Άλγεβρας, οι μετασχηματισμοί αυτοί περιγράφονται με τη βοήθεια πινάκων ως:

$$\begin{pmatrix} x' \\ y' \end{pmatrix} = \begin{pmatrix} 1 & 0 \\ -u & 1 \end{pmatrix} \begin{pmatrix} x \\ y \end{pmatrix}.$$

Οι παραπάνω μετασχηματισμοί απεικονίζουν ευθείες σε ευθείες, παράλληλες ευθείες σε παράλληλες ευθείες, συγγραμμικά ευθύγραμμα τμήματα $AB, \Gamma\Delta$ σε συγγραμμικά ευθύγραμμα τμήματα $A'B', \Gamma'\Delta'$ ώστε $\frac{AB}{\Gamma\Delta} = \frac{A'B'}{\Gamma'\Delta'}$ και σχήματα σε σχήματα με το ίδιο εμβαδόν. Η απόσταση $d_{AB}$ μεταξύ δύο σημείων $A(x_1, y_1)$ και $B(x_2, y_2)$ ορίζεται από τον τύπο

$$d_{AB} = x_2 - x_1$$

είναι δηλαδή το σημασμένο μήκος της προβολής του ευθυγράμμου τμήματος $AB$ επάνω στον άξονα των $x$. Προφανώς η απόσταση μπορεί να πάρει και αρνητικές τιμές, είναι δε

$$d_{BA} = -d_{AB}.$$

Αν η απόσταση $d_{AB}$ μεταξύ των σημείων $A(x_1, y_1)$ και $B(x_2, y_2)$ είναι μηδέν, δηλαδή αν $x_1 = x_2$, τότε τα $A$ και $B$ ανήκουν σε μια κατακόρυφη ευθεία και οι προβολές τους ταυτίζονται. Σε αυτή την περίπτωση η (ειδική) απόστασή τους δίνεται από τον τύπο

$$\delta_{AB} = y_2 - y_1.$$

Επομένως στη Γεωμετρία Galileo δύο σημεία ταυτίζονται μόνο όταν τόσο η απόστασή τους όσο και η ειδική απόστασή τους είναι ίσες με μηδέν.

Για τον ορισμό της γωνίας δύο ευθειών, είναι απαραίτητος ο ορισμός του κύκλου. Έτσι λοιπόν, *κύκλος* ονομάζεται το σύνολο των σημείων $M(x, y)$ του επιπέδου των οποίων οι αποστάσεις από ένα σταθερό σημείο έχουν σταθερή απόλυτη τιμή. Ο ορισμός αυτός οδηγεί στην απεικόνιση του κύκλου στο επίπεδο ως δύο παράλληλες κατακόρυφες ευθείες, ενώ τα (άπειρα) κέντρα του κατοικούν στη μεσοπαράλληλο των δύο ευθειών.

Κατά συνέπεια, για να ορίσουμε τη γωνία δύο ευθειών $ε_1$ και $ε_2$ θεωρούμε τον μοναδιαίο κύκλο, με Κ ένα από τα κέντρα του. Αν θεωρήσουμε ότι οι ευθείες $ε_1$ και $ε_2$ τέμνονται στο Κ, τότε η γωνία τους ορίζεται ως το μήκος του κυκλικού τόξου $N_1N_2$ επί του μοναδιαίου κύκλου. Με την έκφραση «μήκος του κυκλικού τόξου $N_1N_2$» εννοούμε την ειδική απόσταση $δ_{N_1N_2}$. Αποδεικνύεται ότι η γωνία $φ$ δύο ευθειών $ε_1$ και $ε_2$ που έχουν συντελεστές διεύθυνσης $λ_1$ και $λ_2$ αντίστοιχα, δίνεται από τον τύπο

$$φ = δ_{ε_1ε_2} = λ_2 - λ_1.$$

Η Γεωμετρία Minkowski (στις δύο διαστάσεις) είναι η μελέτη των ιδιοτήτων του επιπέδου που παραμένουν αναλλοίωτες κάτω από τη δράση των μετασχηματισμών Lorentz. Ένας τρόπος για την παραγωγή των μετασχηματισμών αυτών, αντίστροφα από την ιστορική ροή αλλά με ξεκάθαρο φυσικό νόημα, είναι να ξεκινήσουμε από τα δύο αξιώματα της Ειδικής Θεωρίας της Σχετικότητας. Θεωρούμε λοιπόν ότι α) οι νόμοι της Φυσικής είναι ίδιοι σε κάθε Α.Σ.Α. και β) η ταχύτητα του φωτός στο κενό ($c$) είναι ίδια σε όλα τα Α.Σ.Α. και ανεξάρτητη από την κίνηση της φωτεινής πηγής. Στη συνέχεια απαιτούμε οι εξισώσεις των μετασχηματισμών να είναι *συμμετρικές* ως προς τα δύο Α.Σ.Α. και *γραμμικές* (Rindler 2003, pp. 11-13). Αποδεικνύεται ότι οι μετασχηματισμοί Lorentz με παράμετρο $u$, όπου $|u| < c = 1$, στο επίπεδο $\{x, t\}$, δηλαδή στον δισδιάστατο χωροχρόνο, είναι:

$$x' = \frac{1}{\sqrt{1-u^2}}x - \frac{u}{\sqrt{1-u^2}}y, \quad y' = -\frac{u}{\sqrt{1-u^2}}x + \frac{1}{\sqrt{1-u^2}}y.$$

Οι παραπάνω μετασχηματισμοί μπορούν να θεωρηθούν *υπερβολική στροφή* (βλ. Callahan 2000, pp. 45-46) και να εκφραστούν με τη βοήθεια πινάκων:

$$\begin{pmatrix} x' \\ y' \end{pmatrix} = \begin{pmatrix} cosh θ & sinh θ \\ sinh θ & cosh θ \end{pmatrix} \begin{pmatrix} x \\ y \end{pmatrix}.$$

Αν ληφθεί υπόψη ότι $cosh^2 θ - sinh^2 θ = 1$, τότε εύκολα αποδεικνύεται ότι οι τελευταίοι μετασχηματισμοί διατηρούν την ποσότητα $x^2 - y^2$ *αναλλοί-* ωτη, δηλαδή αν $x^2 - y^2 = m$ και $\begin{pmatrix} x' \\ y' \end{pmatrix} = \begin{pmatrix} cosh\theta & sinh\theta \\ sinh\theta & cosh\theta \end{pmatrix} \begin{pmatrix} x \\ y \end{pmatrix}$, τότε και $(x')^2 - (y')^2 = m$. Επομένως η ποσότητα $x^2 - y^2$ θα μπορούσε να αντιστοιχεί ή τουλάχιστον να αποτελεί μέρος της νόρμας. Το γεγονός εντούτοις της μη διατήρησης σταθερού προσήμου, εγείρει δυσκολίες. Κατά συνέπεια τα διανύσματα πρέπει να κατανεμηθούν σε κατηγορίες (είδη), ανάλογα με το χωρίο του επιπέδου στο οποίο κατοικούν. Συγκεκριμένα, ένα μη μηδενικό διάνυσμα $a = (x, y)$ καλείται: α) *χωροειδές* αν $x^2 - y^2 > 0$ β) *χρονοειδές* αν $x^2 - y^2 < 0$ και γ) *φωτοειδές* αν $x^2 - y^2 = 0$. Το σύνολο όλων των φωτοειδών διανυσμάτων καλείται *κώνος φωτός*. Όπως και στο Ευκλείδειο επίπεδο, δύο διανύσματα $a, b$ λέγονται *κάθετα* αν $a \cdot b = 0$. Επομένως ο κώνος φωτός αποτελείται από όλα τα διανύσματα που είναι *κάθετα στον εαυτό τους*.

Μπορούμε τώρα να ορίσουμε το *εσωτερικό γινόμενο* δύο διανυσμάτων $a = (x_1, y_1)$ και $b = (x_2, y_2)$ ως $a \cdot b = x_1 x_2 - y_1 y_2$ και με τη βοήθεια αυτού, την αντίστοιχη *νόρμα* ενός διανύσματος $a = (x, y)$ ως:

$$\|a\| = \begin{cases} \sqrt{x^2 - y^2}, \alpha \nu \ a \ \chi\omega\rho o\varepsilon \iota \delta \acute{\varepsilon}\varsigma \\ \sqrt{y^2 - x^2}, \alpha \nu \ a \ \chi\rho o\nu o\varepsilon \iota \delta \acute{\varepsilon}\varsigma \\ \quad 0 \quad , \alpha \nu \ a \ \varphi\omega\tau o\varepsilon \iota \delta \acute{\varepsilon}\varsigma \end{cases}$$

Συνεπώς η *απόσταση* δυο διανυσμάτων $a = (x_1, y_1)$ και $b = (x_2, y_2)$ (του ιδίου είδους) στο επίπεδο Minkowski είναι

$$\|a - b\| = \sqrt{|(x_1 - x_2)^2 - (y_1 - y_2)^2|}$$

ενώ η γωνία $\varphi$ που σχηματίζουν, δίνεται από τον τύπο

$$\varphi = arccosh \frac{a \cdot b}{\|a\|\|b\|}.$$

Οι παραπάνω τρεις μετρικές των μηκών θα μπορούσαν να ενοποιηθούν, εφόσον ενοποιηθούν τα εσωτερικά γινόμενα με τη βοήθεια του τύπου

$$<,>: \mathbb{R}^2 \times \mathbb{R}^2 \to \mathbb{R}: (a, b) \mapsto <a, b> := x_1 x_2 + \varepsilon y_1 y_2$$

όπου η παράμετρος ε παίρνει τις τιμές −1, 0, 1 και οδηγεί στο εσωτερικό γινόμενο, και κατ' επέκταση στη μετρική, Minkowski, Galileo και Ευκλείδη αντίστοιχα.

**4. Διδακτικές προεκτάσεις**

Η παιδαγωγική αξία της διδασκαλίας του μαθήματος της Γεωμετρίας στη Δευτεροβάθμια Εκπαίδευση είναι αδιαμφισβήτητη, κυρίως διότι βοηθάει τους μαθητές στην ανάπτυξη της ικανότητας αντίληψης του χώρου και συνδέει άμεσα τα Μαθηματικά με τον πραγματικό κόσμο (Τουμάσης 2004, σ. 335). Παράλληλα, ως σχολικό μάθημα, εκφράζει ένα περιβάλλον μέσα στο οποίο ευδοκιμεί η μύηση στην αποδεικτική διαδικασία και η κατανόηση του τρόπου θεμελίωσης και ανάπτυξης ενός βασικού κλάδου των Μαθηματικών (Θωμαΐδης 1992). Παρά ταύτα, η σχολική εμπειρία δείχνει ότι η διδασκαλία της Γεωμετρίας, κυρίως στο Λύκειο, παρουσιάζει σημαντικές δυσκολίες, με πιο χαρακτηριστικές την κατανόηση των αποδείξεων (Williams 1979; Driscoll 1982) και τη μηχανιστική εκτέλεσή τους (Patronis & Thomaidis 1997).

Ιδιαίτερα στη χώρα μας, η σταδιακή συρρίκνωση του μαθήματος της Ευκλείδειας Γεωμετρίας και η στατική προσέγγιση στην ανάπτυξη του εναπομείναντος περιεχομένου του, αποδυναμώνουν τα επιχειρήματα αμφισβήτησης της ακραίας θέσης του Jean Dieudonné, ο οποίος πριν 62 χρόνια με στόχο την αναμόρφωση των Αναλυτικών Προγραμμάτων για τα Μαθηματικά, διακήρυττε: «A bas Euclide!» («Κάτω ο Ευκλείδης!»). Στην πράξη η Ευκλείδεια Γεωμετρία διδάσκεται αποκομμένη από το φυσικό της νόημα και εν πολλοίς ασυσχέτιστη με τις εμπειρίες των μαθητών. Αυτό έχει ως συνέπεια να δίδεται η εντύπωση ότι συνιστά ένα κλειστό λογικοπαραγωγικό σύστημα, όπως ένα παιχνίδι με ιδιότυπους κανόνες (λ.χ. το σκάκι). Την ίδια στιγμή, ο φορμαλιστικός χαρακτήρας των σχολικών Μαθηματικών σε συνδυασμό με την πραγματιστική αντίληψη που διαπνέει πολλά εκπαιδευτικά βιβλία και άρθρα, ακόμα και τα λεγόμενα "εκλαϊκευτικά", δεν ευνοούν την ανάπτυξη απόψεων διαφορετικών για την Ευκλείδεια Γεωμετρία πέραν αυτής που είδαμε στην προηγούμενη πρόταση. Έτσι, υιοθετείται εν τοις πράγμασι ένα φορμαλιστικό μοντέλο διδασκαλίας και μάθησης, το οποίο προτάσσει τον δασκαλοκεντρισμό έναντι της ανακάλυψης και ανταλλάσει την εννοιολογική κατανόηση με τη μηχανιστική επίλυση ασκήσεων και τον αλγοριθμικό τρόπο σκέψης.

Από την άλλη μεριά, το πρόγραμμα σπουδών και οι αναλυτικές οδηγίες – ιδιαίτερα την περίοδο της πανδημίας Covid-19 και της συνακόλουθης εφαρμογής της εξ αποστάσεως εκπαίδευσης, οπότε και άλλαξε άρδην το διδακτικό μοντέλο– θέτουν αυστηρούς χρονικούς περιορισμούς στους διδά-σκοντες, οι οποίοι καθημερινά προσπαθούν να ισορροπήσουν μεταξύ κάλυψης της διδακτέας ύλης και κατανόησης από τη μεριά των μαθητών. Σε ένα τέτοιο, σύνθετο, περιβάλλον λοιπόν, το ενδεχόμενο να αφιερωθεί διδακτικός χρόνος στα είδη των Γεωμετριών και στις μεθόδους θεμελίωσής τους πιθανόν να μην φαντάζει ρεαλιστικό.

Έρευνες ωστόσο έχουν δείξει (Cobb 1986; Lenart 1993; van den Brink 1993, 1994; Patronis 1994; Kaisari & Patronis 2010; Rizos et al. 2017) ότι ερχόμενοι οι μαθητές με κατάλληλα μέσα (π.χ. διδακτικά σενάρια ή νοητικά πειράματα) σε επαφή με γεωμετρίες και πέραν της Ευκλείδειας και ανακαλύπτοντας την έννοια και τον ρόλο της μετρικής, έχουν την ευκαιρία να αμφισβητήσουν παγιωμένες αντιλήψεις για τη φύση της ίδιας της Γεωμετρίας, να αλληλεπιδράσουν, να διαπραγματευτούν μαθηματικά νοήματα, να ανακαλύψουν συνδέσεις με τη Φυσική και να κάνουν υπερβάσεις, διαμορφώνοντας έτσι κατάλληλες συνθήκες και ανοίγοντας ορίζοντες για περαιτέρω μάθηση και κατανόηση στα Μαθηματικά. Στο πλαίσιο αυτό είναι ενδεχομένως εφικτή και η ανάγνωση των τριών γεωμετριών που παρουσιάστηκαν στην προηγούμενη παράγραφο, καθώς και η ενσωμάτωση στοιχείων τους στη διδασκαλία. Έτσι μπορεί να δοθεί στους μαθητές ή/και πρωτοετείς φοιτητές η μελλοντική προοπτική της μελέτης α) του τρόπου θεώρησης και του ρόλου των αξιωμάτων, όπως αναπτύχθηκε στην πρώτη παράγραφο β) της κεντρικής έννοιας του «χώρου» και της μοντελοποίησής του γ) της αλληλεπίδρασης της μαθηματικής θέασης με τη φυσική θεωρία δ) της σημασίας των μετασχηματισμών και των αναλλοιώτων τους, άρα και των διαφορετικών αναπαραστάσεων ενός σχήματος (π.χ. του "κύκλου") από γεωμετρία σε γεωμετρία και ε) επεισοδίων από την Ιστορία των Μαθηματικών.

Τα παραπάνω μπορούν να συμβάλουν στη μετάθεση του ορίζοντα της σκέψης των παιδιών και, μέχρι ενός ορισμένου σημείου, στην αλλαγή του τρόπου που αντιλαμβάνονται τα μαθηματικά πρότυπα και τη λειτουργία του σύμπαντος. Οι μαθητές αποδεσμεύονται από την υποχρέωση να δέχονται ως αξίωμα κάτι το οποίο μπορεί να προκύψει από τη μελέτη των φυσικών φαινομένων. Άλλωστε, όπως μας υπενθυμίζει η Ιστορία της Επιστήμης, η προσέγγιση της Γεωμετρίας μέσα από το πείραμα και την παρατήρηση, σε συνδυασμό με τον μακροχρόνιο φιλοσοφικό και επιστημολογικό μετασχηματισμό της έννοιας του «χώρου», οδήγησαν σε επιστημονικές επαναστάσεις και στην υπέρβαση εδραιωμένων φιλοσοφικών πεποιθήσεων για τον φυσικό κόσμο. Θα είχε λοιπόν ενδιαφέρον, μελλοντικά, η διεξαγωγή μιας έρευνας που θα εξέταζε το εάν και με ποιους τρόπους η διδασκαλία της ιστορικής εξέλιξης γεωμετρικών εννοιών και διαδικασιών (όπως αυτές που είδαμε στις παραγράφους 1 και 2) μπορεί να δημιουργήσει κατάλληλες συνθήκες, προκειμένου οι μαθητές Λυκείου να αποκτήσουν θετική στάση για τη Γεωμετρία και να εμπλακούν ενεργά στη διαπραγμάτευση μαθηματικών νοημάτων.

**Βιβλιογραφικές αναφορές**